\documentclass[12pt]{article}
\usepackage{amsfonts,amssymb,subeqnarray,eqsection,indent,epsf}

\begin{document}

\title{\vspace*{-2cm}
       The Brown--Colbourn Conjecture \\
       on Zeros of Reliability Polynomials \\
       is False}

\author{
  \\
  {\small Gordon Royle}                                    \\[-2mm]
  {\small\it Department of Computer Science \& Software Engineering}  \\[-2mm]
  {\small\it University of Western Australia} \\[-2mm]
  {\small\it 35 Stirling Highway} \\[-2mm]
  {\small\it Crawley, WA 6009, AUSTRALIA}                         \\[-2mm]
  {\small\tt GORDON@CSSE.UWA.EDU.AU}   \\[5mm]
  {\small Alan D.~Sokal}                  \\[-2mm]
  {\small\it Department of Physics}       \\[-2mm]
  {\small\it New York University}         \\[-2mm]
  {\small\it 4 Washington Place}          \\[-2mm]
  {\small\it New York, NY 10003 USA}      \\[-2mm]
  {\small\tt SOKAL@NYU.EDU}               \\[-2mm]
  {\protect\makebox[5in]{\quad}}  
  \\
}

\date{January 18, 2003 \\[1mm] revised April 2, 2004}
\maketitle
\thispagestyle{empty}   

\begin{abstract}
We give counterexamples to the Brown--Colbourn conjecture
on reliability polynomials,
in both its univariate and multivariate forms.
The multivariate Brown--Colbourn conjecture is false
already for the complete graph $K_4$.
The univariate Brown--Colbourn conjecture is false
for certain simple planar graphs obtained from $K_4$
by parallel and series expansion of edges.
We show, in fact,
that a graph has the multivariate Brown--Colbourn property
if and only if it is series-parallel.
\end{abstract}


\bigskip
\noindent
{\bf Mathematics Subject Classification (MSC 2000) codes:}
05C99 (Primary);
05C40, 68M10, 68M15, 68R10, 82B20, 90B15, 90B18, 90B25, 94C15 (Secondary).

\bigskip
\noindent
{\bf Running Head:}  Brown--Colbourn Conjecture is False

\clearpage

\newtheorem{defin}{Definition}[section]
\newtheorem{definition}[defin]{Definition}
\newtheorem{prop}[defin]{Proposition}
\newtheorem{proposition}[defin]{Proposition}
\newtheorem{lem}[defin]{Lemma}
\newtheorem{lemma}[defin]{Lemma}
\newtheorem{guess}[defin]{Conjecture}
\newtheorem{ques}[defin]{Question}
\newtheorem{question}[defin]{Question}
\newtheorem{prob}[defin]{Problem}
\newtheorem{problem}[defin]{Problem}
\newtheorem{thm}[defin]{Theorem}
\newtheorem{theorem}[defin]{Theorem}
\newtheorem{cor}[defin]{Corollary}
\newtheorem{corollary}[defin]{Corollary}
\newtheorem{conj}[defin]{Conjecture}
\newtheorem{conjecture}[defin]{Conjecture}

\newtheorem{pro}{Problem}
\newtheorem{clm}{Claim}
\newtheorem{con}{Conjecture}

%
%
\newcounter{example}[section]
\newenvironment{example}%
{\refstepcounter{example}
 \bigskip\par\noindent{\bf Example \thesection.\arabic{example}.}\quad
}%
{\quad $\Box$}
\def\bexam{\begin{example}}
\def\eexam{\end{example}}

\renewcommand{\theenumi}{\alph{enumi}}
\renewcommand{\labelenumi}{(\theenumi)}
\def\prf{\par\noindent{\bf Proof.\enspace}\rm}
\def\rmk{\par\medskip\noindent{\bf Remark.\enspace}\rm}

\newcommand{\be}{\begin{equation}}
\newcommand{\ee}{\end{equation}}
\newcommand{\<}{\langle}
\renewcommand{\>}{\rangle}
\newcommand{\widebar}{\overline}
\def\reff#1{(\protect\ref{#1})}
\def\spose#1{\hbox to 0pt{#1\hss}}
\def\ltapprox{\mathrel{\spose{\lower 3pt\hbox{$\mathchar"218$}}
 \raise 2.0pt\hbox{$\mathchar"13C$}}}
\def\gtapprox{\mathrel{\spose{\lower 3pt\hbox{$\mathchar"218$}}
 \raise 2.0pt\hbox{$\mathchar"13E$}}}
\def\textprime{${}^\prime$}
\def\proof{\par\medskip\noindent{\sc Proof.\ }}
\newcommand{\qed}{\quad $\Box$ \medskip \medskip}
\def\proofof#1{\bigskip\noindent{\sc Proof of #1.\ }}
\def\half{ {1 \over 2} }
\def\third{ {1 \over 3} }
\def\twothird{ {2 \over 3} }
\def\smfrac#1#2{\textstyle{#1\over #2}}
\def\smhalf{ \smfrac{1}{2} }
\newcommand{\real}{\mathop{\rm Re}\nolimits}
\renewcommand{\Re}{\mathop{\rm Re}\nolimits}
\newcommand{\imag}{\mathop{\rm Im}\nolimits}
\renewcommand{\Im}{\mathop{\rm Im}\nolimits}
\newcommand{\sgn}{\mathop{\rm sgn}\nolimits}
\def\hboxscript#1{ {\hbox{\scriptsize\em #1}} }

\newcommand{\restrict}{\upharpoonright}
\renewcommand{\emptyset}{\varnothing}

\def\Z{{\mathbb Z}}
\def\ZZ{{\mathbb Z}}
\def\R{{\mathbb R}}
\def\C{{\mathbb C}}
\def\CC{{\mathbb C}}
\def\N{{\mathbb N}}
\def\NN{{\mathbb N}}
\def\Q{{\mathbb Q}}

\newcommand{\scra}{{\mathcal{A}}}
\newcommand{\scrb}{{\mathcal{B}}}
\newcommand{\scrc}{{\mathcal{C}}}
\newcommand{\scrf}{{\mathcal{F}}}
\newcommand{\scrg}{{\mathcal{G}}}
\newcommand{\scrh}{{\mathcal{H}}}
\newcommand{\scrl}{{\mathcal{L}}}
\newcommand{\scro}{{\mathcal{O}}}
\newcommand{\scrp}{{\mathcal{P}}}
\newcommand{\scrr}{{\mathcal{R}}}
\newcommand{\scrs}{{\mathcal{S}}}
\newcommand{\scrt}{{\mathcal{T}}}
\newcommand{\scrv}{{\mathcal{V}}}
\newcommand{\scrw}{{\mathcal{W}}}
\newcommand{\scrz}{{\mathcal{Z}}}
\newcommand{\scrbt}{{\mathcal{BT}}}
\newcommand{\scrbf}{{\mathcal{BF}}}


\newenvironment{sarray}{
	  \textfont0=\scriptfont0
	  \scriptfont0=\scriptscriptfont0
	  \textfont1=\scriptfont1
	  \scriptfont1=\scriptscriptfont1
	  \textfont2=\scriptfont2
	  \scriptfont2=\scriptscriptfont2
	  \textfont3=\scriptfont3
	  \scriptfont3=\scriptscriptfont3
	\renewcommand{\arraystretch}{0.7}
	\begin{array}{l}}{\end{array}}

\newenvironment{scarray}{
	  \textfont0=\scriptfont0
	  \scriptfont0=\scriptscriptfont0
	  \textfont1=\scriptfont1
	  \scriptfont1=\scriptscriptfont1
	  \textfont2=\scriptfont2
	  \scriptfont2=\scriptscriptfont2
	  \textfont3=\scriptfont3
	  \scriptfont3=\scriptscriptfont3
	\renewcommand{\arraystretch}{0.7}
	\begin{array}{c}}{\end{array}}

\section{Introduction}

Let us consider a connected (multi)graph\footnote{
    Henceforth we omit the prefix ``multi''.
    In this paper a ``graph'' is allowed to have
    loops and/or multiple edges unless explicitly stated otherwise.
}
$G=(V,E)$ as a communications network
with unreliable communication channels,
in which edge $e$ is operational with probability $p_e$
and failed with probability $1-p_e$,
independently for each edge.
Let $R_G({\bf p})$ be the probability that 
every node is capable of communicating with every other node
(this is the so-called {\em all-terminal reliability}\/).
Clearly we have
\be
   R_G({\bf p})  \;=\;
   \sum_{\begin{scarray}
           A \subseteq E \\
           (V,A) \, {\rm connected}
         \end{scarray}}
   \prod_{e \in A} p_e  \prod_{e \in E \setminus A} (1-p_e)
   \;,
\ee
where the sum runs over all connected spanning subgraphs of $G$,
and we have written ${\bf p} = \{p_e\}_{e \in E}$.
We call $R_G({\bf p})$ the (multivariate) {\em reliability polynomial}\/
\cite{Colbourn_87} for the graph $G$;
it is a multiaffine polynomial,
i.e.\ of degree at most 1 in each variable separately.
If the edge probabilities $p_e$ are all set to the same value $p$,
we write the corresponding univariate polynomial as $R_G(p)$,
and call it the univariate reliability polynomial.
We are interested in studying the zeros of these polynomials
when the variables $p_e$ (or $p$) are taken to be {\em complex}\/ numbers.

Brown and Colbourn \cite{Brown_92}
studied a number of examples and made the following conjecture:
\begin{quote}
   {\bf Univariate Brown--Colbourn conjecture.}
   For any graph $G$,
   the zeros of the univariate reliability polynomial $R_G(p)$
   all lie in the closed disc $|p-1| \le 1$.
   In other words, if $|p-1| > 1$, then $R_G(p) \neq 0$.
\end{quote}
Subsequently, one of us \cite{Sokal_chromatic_bounds}
proposed a multivariate extension of the Brown--Colbourn conjecture:
\begin{quote}
   {\bf Multivariate Brown--Colbourn conjecture.}
   For any graph $G$,
   if $\mbox{$|p_e-1| > 1$}$ for all edges $e$, then $R_G({\bf p}) \neq 0$.
\end{quote}

Not long ago, Wagner \cite{Wagner_00} proved,
using an ingenious and complicated construction,
that the univariate Brown--Colbourn conjecture
holds for all series-parallel graphs.\footnote{
   Unfortunately, there seems to be no completely standard
   definition of ``series-parallel graph'';
   a plethora of slightly different definitions can be found in the literature
   \cite{Duffin_65,Colbourn_87,Oxley_86,Oxley_92,Brandstadt_99}.
   So let us be completely precise about our own usage:
   we shall call a loopless graph {\em series-parallel}\/
   if it can be obtained from a forest by a finite sequence of
   series and parallel extensions of edges
   (i.e.\ replacing an edge by two edges in series or two edges in parallel).
   We shall call a general graph (allowing loops) series-parallel
   if its underlying loopless graph is series-parallel.
   Some authors write ``obtained from a tree'', ``obtained from $K_2$''
   or ``obtained from $C_2$'' in place of ``obtained from a forest'';
   in our terminology these definitions yield, respectively,
   all {\em connected}\/ series-parallel graphs,
   all connected series-parallel graphs whose blocks form a path,
   or all {\em 2-connected}\/ series-parallel graphs.
   See \cite[Section 11.2]{Brandstadt_99} for a more extensive bibliography.
}
Subsequently, one of us
\cite[Remark 3 in Section 4.1]{Sokal_chromatic_bounds} showed,
by a two-line induction,
that the multivariate Brown--Colbourn conjecture
holds for all series-parallel graphs.\footnote{
   This proof is reproduced here as
   Theorem~\ref{thm5.6}(c) $\Longrightarrow$ (a).
}
Both the univariate and multivariate conjectures remained open
for general graphs, but most workers in the field
suspected that they would be true.
(At least the present authors did.)

In this short note we would like to report that both the
univariate and multivariate Brown--Colbourn conjectures are false!
The multivariate conjecture is false already for the simplest
non-series-parallel graph, namely the complete graph $K_4$.
As a corollary we will deduce that the univariate conjecture
is false for a 4-vertex, 16-edge planar graph that can be obtained from $K_4$
by adding parallel edges, and for a 1512-vertex, 3016-edge simple planar graph
that can be obtained from $K_4$ by adding parallel edges and then
subdividing edges.
So the Brown--Colbourn conjecture is not true even for simple planar graphs.

Furthermore, for the multivariate property we are able to obtain
a complete characterization:
a graph has the multivariate Brown--Colbourn property
{\em if and only if}\/ it is series-parallel.

It is convenient to restate the Brown--Colbourn conjectures
in terms of the generating polynomial for connected spanning subgraphs,
\be
   C_G({\bf v})  \;=\;
   \sum_{\begin{scarray}
           A \subseteq E \\
           (V,A) \, {\rm connected}
         \end{scarray}}
   \prod_{e \in A} v_e
   \;,
\ee
where we have written ${\bf v} = \{v_e\}_{e \in E}$.
This is clearly related to the reliability polynomial by
\begin{eqnarray}
   R_G({\bf p})  & = &   \left[ \prod_{e \in E} (1-p_e) \right]
                      C_G \!\left( { {\bf p} \over {\bf 1} - {\bf p}} \right)
        \\[4mm]
   C_G({\bf v})  & = &   \left[ \prod_{e \in E} (1+v_e) \right] 
                      R_G \!\left( { {\bf v} \over {\bf 1} + {\bf v}} \right)
\end{eqnarray}
where ${\bf 1}$ denotes the vector with all entries 1,
and division of vectors is understood componentwise.
The multivariate Brown--Colbourn conjecture then states that
if $G$ is a {\em loopless}\/ graph
and $|1+v_e| < 1$ for all edges $e$, then $C_G({\bf v}) \neq 0$.
Loops must be excluded because a loop $e$ multiplies $C_G$
by a factor $1+v_e$ but leaves $R_G$ unaffected.
Some workers also prefer to use the failure probabilities $q_e = 1-p_e$
as the variables.


The plan of this paper is as follows:
In Section~\ref{sec2} we show that the multivariate Brown--Colbourn conjecture
fails for the complete graph $K_4$.
In Section~\ref{sec3} we review the series and parallel reduction formulae
for the reliability polynomial.
In Section~\ref{sec4} we show that the univariate Brown--Colbourn conjecture
fails for certain graphs that are obtained from $K_4$ by adding parallel edges
and then optionally subdividing edges.
In Section~\ref{sec5}
we complete these results by showing
that a graph has the multivariate Brown--Colbourn property
if and only if it is series-parallel.

\section{The multivariate Brown--Colbourn conjecture is false for $K_4$}
   \label{sec2}

For the complete graph $K_4$, the univariate polynomial $C_G(v)$ is
\be
   C_{K_4}(v)  \;=\;  16v^3 + 15v^4 + 6v^5 + v^6  \;.
\ee
The roots of this polynomial all lie outside the disc $|1+v| < 1$,
so the univariate Brown--Colbourn conjecture is true for $K_4$.

\begin{figure}[p]
\vspace*{-2.5cm} \hspace*{0cm}
\begin{center}
\epsfxsize=0.61\textwidth
\epsffile{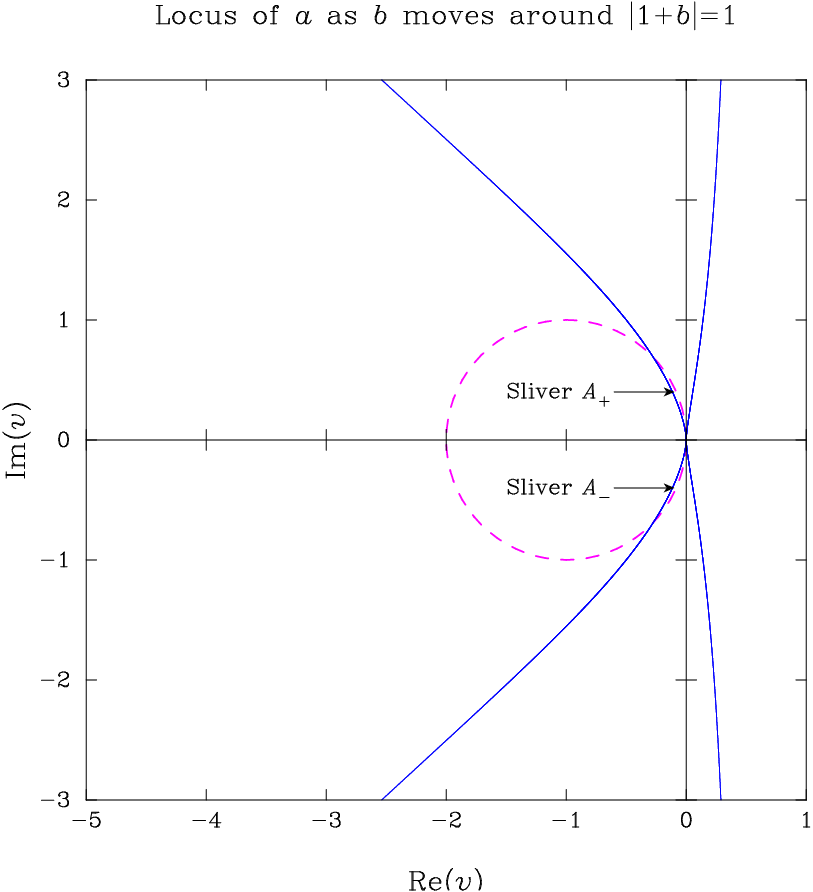} \\
\vspace{1cm}
\epsfxsize=0.61\textwidth
\epsffile{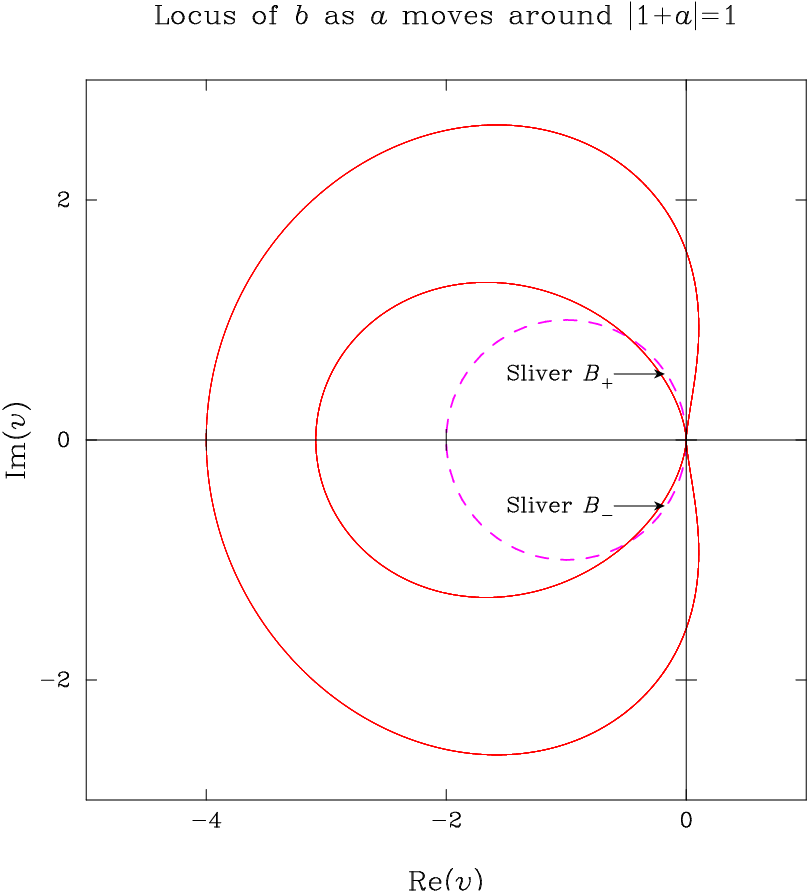} \\
\vspace{-2mm}
\end{center}
\caption{
   Curves for case (b).
   First plot shows the $a$-plane;
   second plot shows the $b$-plane.
   Dashed magenta curve is the circle $|1+v|=1$;
   solid blue curve is the locus of root $a$;
   solid red curve is the locus of root $b$.
}
\label{fig1}
\end{figure}

\begin{figure}[p]
\vspace*{-1.5cm} \hspace*{0cm}
\begin{center}
\leavevmode\epsffile{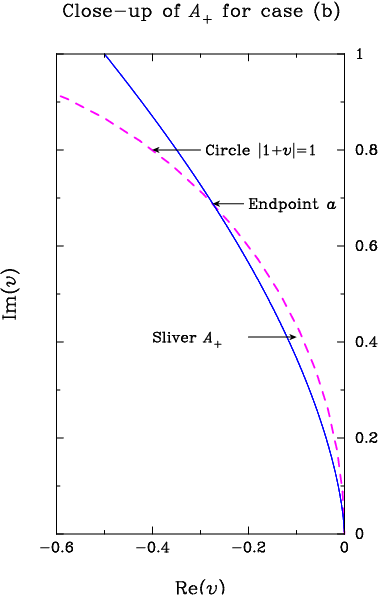}
\hspace{1.5cm}
\epsffile{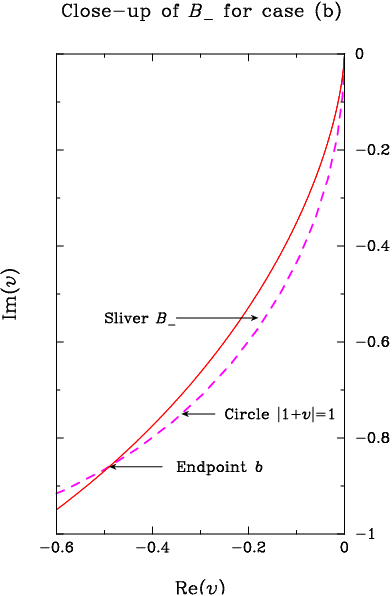} \\
\vspace{1cm}
\leavevmode\epsffile{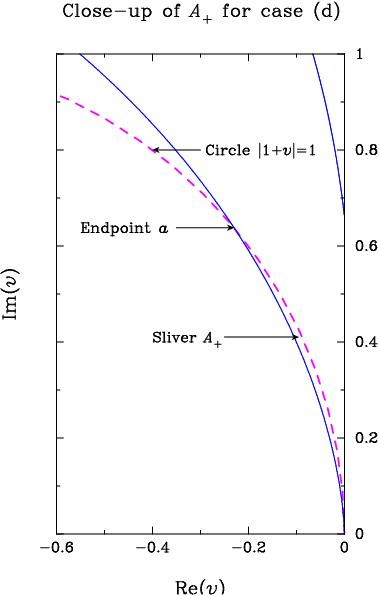}
\hspace{1.5cm}
\epsffile{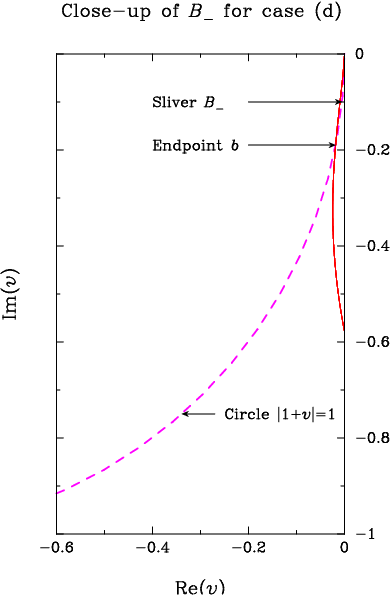} \\
\end{center}
\caption{
   Blow-up of curves to show more clearly the ``sliver'' regions
   $A_+$ and $B_-$.
   Top row shows the $a$- and $b$-planes for case (b);
   bottom row shows the $a$- and $b$-planes for case (d).
}
\label{fig2}
\end{figure}

Let us now consider the bivariate situation,
in which the six edges receive two different weights $a$ and $b$.
There are five cases:
\begin{itemize}
   \item[(a)]  One edge receives weight $a$ and the other five
      receive weight $b$:
\be
   C_{K_4}(a,b)  \;=\;  (8b^3 + 5b^4 + b^5) + (8b^2 + 10b^3 + 5b^4 + b^5) a
\ee
   \item[(b)]  A pair of nonintersecting edges receive weight $a$
      and the other four edges receive weight $b$:
\be
   C_{K_4}(a,b)  \;=\;  (4b^3 + b^4) + (8b^2 + 8 b^3 + 2b^4) a
                           + (4b + 6b^2 + 4 b^3 + b^4) a^2
\ee
   \item[(c)]  A pair of intersecting edges receive weight $a$
      and the other four edges receive weight $b$:
\be
   C_{K_4}(a,b)  \;=\;  (3b^3 + b^4) + (10b^2 + 8 b^3 + 2b^4) a
                           + (3b + 6b^2 + 4 b^3 + b^4) a^2
\ee
   \item[(d)]  A 3-star receives weight $a$ and the complementary
      triangle receives weight $b$:
\be
   C_{K_4}(a,b)  \;=\;  (9b^2 + 3b^3) a + (6b + 9b^2 + 3b^3) a^2
                           + (1 + 3b + 3b^2 + b^3) a^3
\ee
   \item[(e)]  A three-edge path receives weight $a$ and the
      complementary three-edge path receives weight $b$:
\be
   C_{K_4}(a,b)  \;=\;  b^3 + (7b^2+3b^3)a + (7b+9b^2+3b^3)a^2
                           + (1+3b+3b^2+b^3)a^3
\ee
\end{itemize}
We have plotted the roots $a$ when $b$ traces out the circle $|1+b|=1$,
and vice versa.
In cases (b) and (d) it turns out that the roots
can enter the ``forbidden discs'' $|1+a| < 1$ and $|1+b| < 1$.
This is shown in Figure~\ref{fig1} for case (b);
blow-ups of the crucial regions are shown in Figure~\ref{fig2}
both for case (b) and for case (d).
As a result, counterexamples to the multivariate Brown--Colbourn conjecture
can be obtained in these two cases:
indeed, for any $a$ lying in the region $A_+$ (resp.\ $A_-$),
there exists $b \in B_-$ (resp.\ $B_+$)
such that $C_{K_4}(a,b) = 0$, and conversely.

Let us note for future reference that
the endpoint of the region $A_\pm$ (resp.\ $B_\pm$)
lies at $a = -1 + e^{\pm 2\pi i \alpha}$
(resp.\ $b = -1 + e^{\pm 2\pi i \beta}$),
where $\alpha \approx 0.120692$ and $\beta \approx 0.164868$ in case (b),
and $\alpha \approx 0.110198$ and $\beta \approx 0.030469$ in case (d).

We can understand this behavior analytically as follows:
For each of the five cases,
let us solve the equation $C_{K_4}(a,b) = 0$ for $a$ in terms of $b$,
expanding in power series for $b$ near 0.
We obtain:
\begin{itemize}
   \item[(a)]  $a = -b + {5 \over 8} b^2 + O(b^3)$
   \item[(b)]  $a = -b \pm {1 \over 2} b^{3/2} + O(b^2)$
   \item[(c)]  $a = -{1 \over 3} b + {1 \over 8} b^2 + O(b^3)$ and
               $a = -3b + {31 \over 8} b^2 + O(b^3)$
   \item[(d)]  $a = -3b \pm i \sqrt{3} b^{3/2} + O(b^2)$ and $a=0$
   \item[(e)]  $a = -b + {3 \over 4} b^2 + O(b^3)$ and
               $a = (-3 \pm 2 \sqrt{2})b + {9 \over 16} (10 \mp 7\sqrt{2})b^2
                    + O(b^3)$
\end{itemize}
The behavior is thus different in cases (a,c,e) on the one hand
and cases (b,d) on the other:

\bigskip

{\bf Cases (a,c,e):}\/
Here the solution is of the form
\be
   a   \;=\;  \gamma_1 b + \gamma_2 b^2 + O(b^3)
\ee
with $\gamma_1,\gamma_2$ real.
Therefore, if we set $b = -1 + e^{i\theta}$ and expand in powers of $\theta$,
we obtain
\be
   |1+a|^2  \;=\;  1 \,+\, (\gamma_1^2 - \gamma_1 - 2\gamma_2) \theta^2
                     \,+\, O(\theta^4)
   \;.
\ee
Provided that $\gamma_1^2 - \gamma_1 - 2\gamma_2 > 0$
--- as indeed holds for all the roots in cases (a,c,e) ---
we have $|1+a| \ge 1$ for small $\theta$,
so no counterexample is found (at least for small $\theta$).

\bigskip

{\bf Cases (b,d):}\/
Here, by contrast, the solution is of the form
\be
   a   \;=\;  \delta_1 b + \delta_2 b^{3/2} + O(b^2)
\ee
with $\delta_1 < 0$ and $\delta_2 \neq 0$.
Therefore, if we set $b = -1 + e^{i\theta}$ and expand as before,
we obtain
\be
   a   \;=\;  i \delta_1 \theta \,+\, e^{\pm 3\pi i/4} \delta_2 \theta^{3/2}
              \,+\, O(\theta^2)
   \;.
\ee
Since $\real(e^{\pm 3\pi i/4} \delta_2) < 0$ for at least one of the roots,
we have $\real a \propto -|\imag a|^{3/2}$ for small $\theta$;
in particular, we have $|1+a| < 1$ for small $\theta \neq 0$.

In fact, more can be said:
suppose that we fix {\em any}\/ $\lambda > 0$
and set $b = \lambda(-1 + e^{i\theta})$.
Then we have
\be
   a  \;=\;  i \delta_1 \lambda \theta \,+\,
             e^{\pm 3\pi i/4} \delta_2 \lambda^{3/2} \theta^{3/2}
             \,+\, O(\theta^2)
   \;,
\ee
so that once again $\real a \propto -|\imag a|^{3/2}$ for small $\theta$.
In particular, we will have $|\lambda+a| < \lambda$ for small $\theta \neq 0$,
irrespective of how small $\lambda$ was chosen.
This observation will play a crucial role in Section~\ref{sec5}
(see Proposition~\ref{prop5.5}).

\section{Series and parallel reduction formulae}
   \label{sec3}

Suppose that $G$ contains edges $e_1,\ldots,e_n$
(with corresponding weights $v_1,\ldots,v_n$)
in parallel between the same pair of vertices $x,y$.
Then it is easy to see that the edges $e_1,\ldots,e_n$
can be replaced by a single edge of weight
\be
   v_1 \| v_2 \| \cdots \| v_n
   \;\equiv\;
   \prod_{i=1}^n (1+v_i) \,-\, 1
 \label{eq.parallel}
\ee
without changing the value of $C_G({\bf v})$.
[Reason: $x$ is connected to $y$ via this ``super-edge''
 if and only if $x$ is connected to $y$ by at least one
 of the edges $e_1,\ldots,e_n$.]

Suppose next that $G$ contains edges $e_1,\ldots,e_n$
(with corresponding weights $v_1,\ldots,v_n$)
in series between the pair of vertices $x,y$:
this means that the edges $e_1,\ldots,e_n$ form a path
in which all the vertices except possibly the
endvertices $x$ and $y$ have degree 2 in $G$.
Let $G'$ be the graph in which the edges $e_1,\ldots,e_n$
are replaced by a single edge $e_*$ from $x$ to $y$.
Then it is not hard to see that
\be
   C_G({\bf v})  \;=\;  \left( \sum_{j=1}^n \prod_{i \neq j} v_i \right)
                        C_{G'}({\bf v'})
 \label{eq.series0}
\ee
where the edge $e_*$ is given weight
\be
   v'_*  \;=\;
   v_1 \bowtie v_2 \bowtie \cdots \bowtie v_n
   \;\equiv\;
   {1 \over \sum\limits_{i=1}^n  1/v_i}
 \label{eq.series}
\ee
and all edges other than $e_1,\ldots,e_n,e_*$ are given weight $v'_e = v_e$.
[Reason:  A connected spanning subgraph of $G$ can omit at most one
of the edges $e_1,\ldots,e_n$, for otherwise at least one of the internal
vertices of the path would be disconnected from both $x$ and $y$.
Moreover, $x$ is connected to $y$ via the ``super-edge'' $e_*$
if and only if none of the edges $e_1,\ldots,e_n$ are omitted.
The relative weight of the cases with and without $x$ connected to $y$
via $e_*$ is thus
$(\prod_{i=1}^n v_i) / (\sum_{j=1}^n \prod_{i \neq j} v_i) = v_*$;
and there is an overall normalization factor
$\sum_{j=1}^n \prod_{i \neq j} v_i$.
See also \cite[p.~35]{Colbourn_87} for an equivalent formula.]

The formula for series reduction can be applied immediately
to handle arbitrary subdivisions of a graph $G$.
Given a finite graph $G=(V,E)$ and a family of integers
${\bf s} = \{s_e\}_{e \in E} \ge 1$,
we define $G^{\bowtie \bf s}$ to be the graph in which each edge $e$ of $G$
is subdivided into $s_e$ edges in series.
If $s \ge 1$ is an integer,
we define $G^{\bowtie s}$ to be the graph in which each edge of $G$
is subdivided into $s$ edges in series.
All the edges in $G^{\bowtie \bf s}$ or $G^{\bowtie s}$
obtained by subdividing the edge $e \in E$
are assigned the same weight $v_e$ as was assigned to $e$ in the
original graph $G$.
It follows immediately from \reff{eq.series0}/\reff{eq.series} that
\be
   C_{G^{\bowtie \bf s}}({\bf v})  \;=\;
   \left( \prod_{e \in E} s_e v_e^{s_e-1} \right)   C_G( {\bf v}/{\bf s} )
 \label{eq.subdivision}
\ee
where $({\bf v}/{\bf s})_e \equiv v_e/s_e$.

\bigskip

{\bf Remarks.}  1.  Series and parallel reduction formulae can be
derived in the more general context of the $q$-state Potts model
(also known as the multivariate Tutte polynomial):
see e.g.\ \cite[Section 2]{Sokal_chromatic_roots}.
Parallel reduction is always given by \reff{eq.parallel},
independently of the value of the parameter $q$.
Series reduction is given by
\be
   v_1 \bowtie v_2 \bowtie \cdots \bowtie v_n
   \;=\;
   {q \over \prod\limits_{i=1}^n (1 + q/v_i) \,-\, 1}
   \;.
 \label{eq.series.q}
\ee
Please note that \reff{eq.series.q} reduces to \reff{eq.series}
when $q \to 0$, which is precisely the limit in which the
multivariate Tutte polynomial $Z_G(q,{\bf v})$ tends
(after division by $q$) to $C_G({\bf v})$.

2.  If one takes in $C_G({\bf v})$ the further limit
of ${\bf v}$ infinitesimal,
one obtains the generating polynomial of
{\em minimal}\/ connected spanning subgraphs,
i.e.\ spanning trees.
Now, spanning trees are intimately related to linear electrical circuits,
as was noticed by Kirchhoff in 1847 \cite{Kirchhoff,Nerode_61}.
For $v$ infinitesimal,
the parallel reduction formula \reff{eq.parallel} becomes
\be
   v_1 \| v_2 \| \cdots \| v_n
   \;\equiv\;
   v_1 + v_2 + \ldots + v_n
   \;,
 \label{eq.parallel.2}
\ee
which is precisely the law for putting electrical conductances in parallel.
And the series reduction formula \reff{eq.series}
is precisely the law for putting electrical conductances in series!

\section{The univariate Brown--Colbourn conjecture is false as well}
   \label{sec4}

Let $K_4^{({\rm a},p_1,p_2)}$ be the graph obtained from $K_4$
by replacing one edge by $p_1$ parallel edges
and replacing each of the other five edges by $p_2$ parallel edges.
Let $K_4^{({\rm b},p_1,p_2)}$ be the graph obtained from $K_4$
by replacing each of two nonintersecting edges by $p_1$ parallel edges
and replacing each of the remaining four edges by $p_2$ parallel edges.
Define in a similar manner
$K_4^{({\rm c},p_1,p_2)}$, $K_4^{({\rm d},p_1,p_2)}$
and $K_4^{({\rm e},p_1,p_2)}$
for the cases (c), (d) and (e) discussed in Section~\ref{sec2}.

We saw in Section~\ref{sec2} that in cases (b) and (d)
one can obtain a counterexample to the multivariate Brown--Colbourn conjecture
by choosing the weight $a$ to lie anywhere in the region $A_+$;
this leads to a root $b$ lying in the region $B_-$
(see Figures~\ref{fig1} and \ref{fig2}).
Note now that the $p$th power of the region $1 + A_+$
will overlap the region $1 + B_-$
whenever $p > (1-\beta)/\alpha$
[just choose any point $b \in B_-$ close enough to the endpoint
 $-1 + e^{-2\pi i \beta} = -1 + e^{2\pi i (1-\beta)}$;
 then one of the $p$th roots of $1+b$ will lie in the region $1 + A_+$].
And \reff{eq.parallel} tells us that $p$ edges in parallel,
each with weight $v$, are equivalent to a single edge with weight
$v_{\rm eff}$ satisfying $1 + v_{\rm eff} = (1+v)^p$.
This reasoning suggests that counterexamples to the univariate
Brown--Colbourn conjecture might be found for the graphs
$K_4^{({\rm b},1,p)}$ and $K_4^{({\rm d},1,p)}$:
for all $p > (1-\beta)/\alpha$ they should have a root $v \in A_+$.\footnote{
   We do not claim that this is a proof,
   though we suspect that a suitable topological argument
   might be able to turn it into a proof.
}
Likewise, the graphs $K_4^{({\rm b},p,1)}$ and $K_4^{({\rm d},p,1)}$
are expected to have, for all $p > (1-\alpha)/\beta$, a root $v \in B_+$.
These guesses are in fact correct, and we find the following
counterexamples to the univariate Brown--Colbourn conjecture:
\begin{itemize}
   \item $G = K_4^{({\rm b},1,7)}$, 30 edges:
             $v \approx -0.269253 \pm 0.682304 i$, $|1 + v| \approx 0.999765$ 
   \item $G = K_4^{({\rm b},6,1)}$, 16 edges:
             $v \approx -0.405015 \pm 0.801589 i$, $|1 + v| \approx 0.998274$
   \item $G = K_4^{({\rm d},1,9)}$, 30 edges:
             $v \approx -0.220759 \pm 0.626655 i$, $|1 + v| \approx 0.999956$
   \item $G = K_4^{({\rm d},30,1)}$, 93 edges:
             $v \approx -0.017476 \pm 0.185846 i$, $|1 + v| \approx 0.999946$
\end{itemize}
Counterexamples are also obtained for each larger $p$;
some typical numbers are shown in Table~\ref{table1}.
Please note that all these counterexample graphs are planar.

\begin{table}[t]
\footnotesize
\hspace*{-4mm}
\begin{tabular}{|c||*{10}{@{\hspace{1mm}}c@{\hspace{1mm}}|}}
\cline{2-11}
\multicolumn{1}{c|}{\quad}  &
 \multicolumn{10}{|c|}{Value of $p$} \\
\hline
Graph  & 6 & 7 & 8 & 9 & 10 & 11 & 12 & 13 & 14 & 15 \\
\hline\hline
$K_4^{({\rm b},1,p)}$ \vphantom{{\Large P}} &
   1 & 0.999765 & 0.997818 & 0.996996 & 0.996734 &
      0.996749 & 0.996897 & 0.997102 & 0.997326 & 0.997547  \\
\hline
$K_4^{({\rm b},p,1)}$ \vphantom{{\Large P}} &
   0.998274 & 0.997234 & 0.997001 & 0.997083 & 0.997284 &
      0.997519 & 0.997753 & 0.997971 & 0.998169 & 0.998345  \\
\hline
$K_4^{({\rm d},1,p)}$ \vphantom{{\Large P}} &
   1 & 1 & 1 & 0.999956 & 0.999813 &
      0.999746 & 0.999718 & 0.999713 & 0.999718 & 0.999730 \\
\hline
$K_4^{({\rm d},p,1)}$ \vphantom{{\Large P}} &
   1 & 1 & 1 & 1 & 1 & 1 & 1 & 1 & 1 & 1 \\
\hline
\end{tabular}
\caption{
   Minimum value of $|1+v|$ for a zero of $C_G(v)$
   for selected graphs $G = K_4^{({\rm b/d},p_1,p_2)}$.
   For $1 \le p \le 5$ the value equals 1.
   A value strictly less than 1 indicates a counterexample
   to the univariate Brown--Colbourn conjecture.
   For $K_4^{({\rm d},p,1)}$ a counterexample can be found for $p \ge 30$.
}
\label{table1}
\end{table}

The graphs $G = K_4^{({\rm b/d},p_1,p_2)}$ are, of course, non-simple
(except when $p_1 = p_2 = 1$);
so one might cling to the hope that the univariate Brown--Colbourn conjecture
is true at least for {\em simple}\/ graphs
(or, weaker yet, for {\em simple planar}\/ graphs).
But these hopes too are false.
To see why this is the case, consider the following procedure:
\begin{itemize}
   \item[1)] Choose $p_1,p_2$ so that the graph $K_4^{({\rm b},p_1,p_2)}$
      has a root $v_1$ satisfying $|1+v_1| < 1$.
   \item[2)] Choose any integer $s \ge 2$.
   \item[3)] Find an integer $k$ large enough so that
        $v_k \equiv -1 + (1+v_1)^{1/k}$ ---
        defined using the root with $|\arg[(1+v_1)^{1/k}]| \le \pi/k$ ---
        lies in the disc $|1/s+v_k| < 1/s$.
        [It is always possible to find such a $k$,
        because the points $v_k$ lie on a logarithmic spiral
        that approaches the point $v=0$
        making a nonzero angle with the imaginary axis,
        while all the circles $|1/s+v|=1/s$
        pass through $v=0$ tangent to the imaginary axis.]
\end{itemize}
Then $v_k$ is a root for the graph $K_4^{({\rm b},kp_1,kp_2)}$,
by the rules for parallel reduction;
and $sv_k$ is a root for the graph $(K_4^{({\rm b},kp_1,kp_2)})^{\bowtie s}$,
by the rules for series reduction.
And by construction we have $|1 + sv_k| < 1$.
Therefore, the graph $(K_4^{({\rm b},kp_1,kp_2)})^{\bowtie s}$,
which is simple and planar,
is the desired counterexample.

%
%
For example, if we take $(p_1,p_2) = (11,1)$ and $s=2$,
counterexamples can be obtained for $k \ge 58$:
\begin{itemize}
   \item $v_1 \approx -0.140\,970\,808\,664 + 0.507\,062\,767\,880 i$,
         $|1+v_1| \approx 0.997\,518\,822\,949$
   \item $v_{58} \approx -0.000\,085\,091\,565 + 0.009\,193\,226\,407 i$,
         $|1 + 2v_{58}| \approx 0.999\,998\,862\,173$
\end{itemize}
This shows that the graph $(K_4^{({\rm b},638,58)})^{\bowtie 2}$,
which has 1512 vertices and 3016 edges,
is a counterexample to the univariate Brown--Colbourn conjecture.
Similarly, if we take $(p_1,p_2) = (1,12)$ and $s=2$,
counterexamples can be obtained for $k \ge 36$:
\begin{itemize}
   \item $v_1 \approx -0.112\,358\,418\,620 + 0.453\,757\,934\,703 i$,
         $|1+v_1| \approx 0.996\,897\,106\,175$
   \item $v_{36} \approx -0.000\,172\,469\,038 + 0.013\,125\,252\,246 i$,
         $|1+2v_{36}| \approx 0.999\,999\,665\,908$
\end{itemize}
Therefore, the graph $(K_4^{({\rm b},36,432)})^{\bowtie 2}$,
which has 1804 vertices and 3600 edges,
is a counterexample to the univariate Brown--Colbourn conjecture.

Smaller counterexamples of the forms
$(K_4^{({\rm b/d},p,1)})^{\bowtie (s,1)}$ or
$(K_4^{({\rm b/d},1,p)})^{\bowtie (1,s)}$
can probably be found by direct search.
But the foregoing construction has the advantage that there is no need
to compute the roots of extremely-high-degree polynomials;
it suffices to compute the roots for the base case $K_4^{({\rm b/d},p_1,p_2)}$
(for which the polynomials are large but not huge)
and then make simple manipulations on them.

\bigskip

{\bf Methodological remark.}
In this work we needed to compute accurately the roots of polynomials
of fairly high degree (up to 93) with very large integer coefficients
(up to about $10^{27}$).
To do this we used the package MPSolve 2.0
developed by Dario Bini and Giuseppe Fiorentino
\cite{Bini_package,Bini-Fiorentino}.
MPSolve is much faster than {\sc Mathematica}'s {\tt NSolve}
for high-degree polynomials
(this is reported in \cite{Bini-Fiorentino}, and we confirm it);
it gives guaranteed error bounds for the roots,
based on rigorous theorems \cite{Bini-Fiorentino};
its algorithms are publicly documented \cite{Bini-Fiorentino};
and its source code is freely available \cite{Bini_package}.

\bigskip

Let us mention, finally, that counterexamples with smaller values of $|1+v|$
can be found.  Consider, for example, the complete graph $K_6$
in which a pair of vertex-disjoint triangles receives weight $a$
and the remaining nine edges receive weight $b$.  We have
\begin{eqnarray}
   &  & C_{K_6}(a,b)  \;=\;
       \nonumber \\
   & & \quad
   (81 b^5 + 78 b^6 + 36 b^7 + 9 b^8 + b^9) +
       \nonumber \\
   & & \quad
   (324 b^4 + 594 b^5 + 480 b^6 + 216 b^7 + 54 b^8 + 6 b^9) a +
       \nonumber \\
   & & \quad
   (486 b^3 + 1314 b^4 + 1665 b^5 + 1224 b^6 + 540 b^7 + 135 b^8 + 15 b^9)
       a^2  + 
       \nonumber \\
   & & \quad
   (324 b^2 + 1188 b^3 + 2160 b^4 + 2376 b^5 + 1656 b^6 + 720 b^7 +
       180 b^8 + 20 b^9) a^3 +
       \nonumber \\
   & & \quad
   (81 b + 432 b^2 + 1134 b^3 + 1800 b^4 + 1854 b^5 +
        1254 b^6 + 540 b^7 + 135 b^8 + 15 b^9) a^4  +
       \nonumber \\
   & & \quad
   (54 b + 216 b^2 + 504 b^3 + 756 b^4 +
      756 b^5 + 504 b^6 + 216 b^7 + 54 b^8 + 6 b^9) a^5  +
       \nonumber \\
   & & \quad
   (9 b + 36 b^2 + 84 b^3 + 126 b^4 + 126 b^5 + 84 b^6 +
      36 b^7 + 9 b^8 + b^9) a^6
\end{eqnarray}
If we then substitute $a = (1+v)^{p_1} - 1$ and $b = (1+v)^{p_2} - 1$,
counterexamples to the univariate Brown--Colbourn conjecture
can be found for many pairs $(p_1,p_2)$.
For example, for $(p_1,p_2) = (1,6)$ we obtain a 60-edge non-planar graph
whose roots include $v \approx -0.357514 \pm 0.713815\,i$,
yielding $|1+v| \approx 0.960375$.

It would be interesting to know whether examples can be found in which
$|1+v|$ is arbitrarily small.  More generally, one can ask:

\begin{question}
What is the closure of the set of all roots of the polynomials $C_G(v)$
as $G$ ranges over all graphs?
Over all planar graphs?
Over all simple planar graphs?
\end{question}

\noindent
Brown and Colbourn \cite{Brown_92} pointed out that the graphs $G = C_n^{(p)}$
(the $n$-cycle with each edge replaced by $p$ parallel edges)
have roots that, taken together, are dense in the region $|1+v| \ge 1$.
We have shown here that roots can also enter the region $|1+v| < 1$.
But how far into this latter region can they penetrate?
Might the roots actually be dense in the whole complex plane?
If this is indeed the case,
it would mean that the univariate Brown--Colbourn conjecture
is as false as it can possibly be.

\bigskip

{\bf Note Added (April 2004).}
Building on the examples constructed in this section,
Chang and Shrock \cite[Sections 5.17 and 5.18]{Shrock_reliability}
have recently devised families of strip graphs
in which the limiting curve of zeros of $C_G(v)$,
as the strip length tends to infinity,
penetrates into the ``forbidden region'' $|1+v| < 1$.
Some of these families consist of planar graphs.

\section{Series-parallel is necessary and sufficient}
   \label{sec5}

In this section we shall prove
that a graph has the multivariate Brown--Colbourn property
{\em if and only if}\/ it is series-parallel.

Let us begin by defining a weakened version of the Brown--Colbourn property:

\begin{definition}
Let $G$ be a graph, and let $\lambda > 0$.
We say that $G$
\begin{itemize}
   \item has the {\em univariate property BC${}_\lambda$}\/
      if $C_G(v) \neq 0$ whenever $|\lambda + v| < \lambda$;
   \item has the {\em multivariate property BC${}_\lambda$}\/
      if $C_G({\bf v}) \neq 0$ whenever $|\lambda + v_e| < \lambda$
      for all edges $e$.
\end{itemize}
\end{definition}
Properties BC${}_1$ are, of course, the original
univariate and multivariate Brown--Colbourn properties;
the properties BC${}_\lambda$ become increasingly weaker
as $\lambda$ is decreased.

The properties BC${}_\lambda$ are intimately related to subdivisions:

\begin{lemma}
   \label{lemma5.2}
Let $\lambda > 0$ and let $s$ be a positive integer.
Then the following are equivalent for a graph $G$:
\begin{itemize}
   \item[(a)]  $G$ has the univariate property BC${}_\lambda$.
   \item[(b)]  $G^{\bowtie s}$ has the univariate property BC${}_{s\lambda}$.
\end{itemize}
\end{lemma}

\begin{lemma}
   \label{lemma5.3}
Let $\lambda > 0$ and let $s$ be a positive integer.
Then the following are equivalent for a graph $G$:
\begin{itemize}
   \item[(a)]  $G$ has the multivariate property BC${}_\lambda$.
   \item[(b)]  $G^{\bowtie s}$ has the multivariate property
       BC${}_{s\lambda}$.
   \item[(c)]  $G^{\bowtie \bf s}$ has the multivariate property
       BC${}_{s\lambda}$ for all vectors ${\bf s}$ satisfying
       $s_e \ge s$ for all edges $e$.
\end{itemize}
\end{lemma}

\noindent
Indeed, Lemmas~\ref{lemma5.2} and \ref{lemma5.3}
are an immediate consequence of the formula \reff{eq.subdivision}
for subdivisions --- which states that subdivision by ${\bf s}$
moves the nonzero roots from ${\bf v}$ to ${\bf sv}$ ---
together with the fact that
$|s\lambda + sv| < s\lambda$ is equivalent to $|\lambda + v| < \lambda$.

In the preceding section we have shown that not all graphs
have the univariate property BC${}_1$.
It is nevertheless true --- and virtually trivial ---
that every connected graph has the univariate property BC${}_\lambda$
for {\em some}\/ $\lambda > 0$.
(Since a non-identically-vanishing univariate polynomial
has finitely many roots,
it suffices to choose $\lambda$ small enough so that
none of the roots of $C_G(v)$ lie in the disc $|\lambda + v| < \lambda$.)
By Lemma~\ref{lemma5.2}, an equivalent assertion is that
$G^{\bowtie s}$ has the univariate property BC${}_1$
for all sufficiently large integers $s$.\footnote{
   Brown and Colbourn \cite[Proposition 4.4 and Theorem 4.5]{Brown_92}
   have proven a result also for \emph{non-uniform} subdivisions
   $G^{\bowtie \bf s}$:
   namely, for each graph $G$ there exists an integer $s$
   such that $G^{\bowtie \bf s}$ has the univariate property BC${}_1$
   whenever $s_e \ge s$ for all $e$.
   This is significantly stronger than the just-mentioned trivial result,
   and it would be worth trying to understand it better.
   Brown and Colbourn's method looks very different from ours,
   at least at first glance;
   it would be interesting to try to translate it into our language.
   In particular, there may be a ``partially multivariate'' result
   hiding underneath their apparently univariate proof.
}

The situation is very different, however,
when we consider the multivariate property BC${}_\lambda$.
We begin with a simple but important lemma:

\begin{lemma}
   \label{lemma5.4}
Let $\lambda > 0$, and suppose that the connected graph $G$
has the multivariate property BC${}_\lambda$.
Then every connected subgraph $H \subseteq G$
also has the multivariate property BC${}_\lambda$.
\end{lemma}

\proof
Consider first the case in which $H$ is a connected {\em spanning}\/ subgraph
(i.e.\ its vertex set is the same as that of $G$).
Let us write ${\bf v} = ({\bf v}', {\bf v}'')$
where ${\bf v}' = \{v_e\}_{e \in E(H)}$
and ${\bf v}'' = \{v_e\}_{e \in E(H) \setminus E(G)}$.
Then
\be
   C_H({\bf v}')  \;=\;  C_G({\bf v}', {\bf 0})  \;=\;
        \lim_{{\bf v}'' \to {\bf 0}}  C_G({\bf v}', {\bf v}'')   \;.
\ee
By hypothesis, $C_G({\bf v}', {\bf v}'') \neq 0$
whenever $|\lambda + v_e| < \lambda$ for all $e \in E(G)$.
Now take ${\bf v''} \to {\bf 0}$ from within this product of discs
(${\bf 0}$ lies on its boundary).
By Hurwitz's theorem\footnote{
   Hurwitz's theorem states that if $D$ is a domain in $\C^n$
   and $(f_k)$ are nonvanishing analytic functions on $D$
   that converge to $f$ uniformly on compact subsets of $D$,
   then $f$ is either nonvanishing or else identically zero.
   Hurwitz's theorem for $n=1$ is proved in most standard texts
   on the theory of analytic functions of a single complex variable
   (see e.g.\ \cite[p.~176]{Ahlfors_66}).
   Surprisingly, we have been unable to find Hurwitz's theorem
   proven for general $n$ in any standard text on several complex variables
   (but see \cite[p.~306]{Krantz_92} and \cite[p.~337]{Simon_74}).
   So here, for completeness, is the sketch of a proof:
   Suppose that $f(c) = 0$ for some $c = (c_1,\ldots,c_n) \in D$,
   and let $D' \subset D$ be a small polydisc centered at $c$.
   Applying the single-variable Hurwitz theorem,
   we conclude that $f(z_1,c_2,\ldots,c_n) = 0$
   for all $z_1$ such that $(z_1,c_2,\ldots,c_n) \in D'$.
   Applying the same argument repeatedly in the variables $z_2,\ldots,z_n$,
   we conclude that $f$ is identically vanishing on $D'$
   and hence, by analytic continuation, also on $D$.
},
either $C_H({\bf v}')$ is nonvanishing whenever
$|\lambda + v_e| < \lambda$ for all $e \in E(H)$,
or else $C_H$ is identically zero.
But the latter is impossible since $H$ is connected.

Now let $H$ be an arbitrary connected subgraph of $G$ (spanning or not).
Construct a connected spanning subgraph $\widehat{H}$ of $G$
by hanging trees off some or all of the vertices of $H$
without creating any new circuits.\footnote{
   This can be done, for instance, by running breadth-first search
   with the vertices of $H$ initially on the queue.
}
Let us write ${\bf v} = \{v_e\}_{e \in E(\widehat{H})} = ({\bf v}', {\bf v}'')$
where ${\bf v}' = \{v_e\}_{e \in E(H)}$
and ${\bf v}'' = \{v_e\}_{e \in E(\widehat{H}) \setminus E(H)}$.
Then
\be
   C_{\widehat{H}}({\bf v})  \;=\;
   C_H({\bf v}') \prod_{e \in E(\widehat{H}) \setminus E(H)} v_e
   \;.
\ee
Since $\widehat{H}$ has multivariate property BC${}_\lambda$,
so does $H$.
\qed

The following is the fundamental fact from which all else flows:

\begin{proposition}
   \label{prop5.5}
The complete graph $K_4$ does not have the multivariate property
BC${}_\lambda$ for any $\lambda > 0$.
\end{proposition}

\proof
This is an almost immediate consequence of the observations
made at the end of Section~\ref{sec2}.
In cases (b) and (d),
for any $\lambda > 0$ there exists $b$ with $|\lambda+b|=\lambda$
for which at least one of the solutions to $C_{K_4}(a,b) = 0$
satisfies $|\lambda+a| < \lambda$.
By slightly perturbing this pair, we can find a pair $(a,b)$
with $C_{K_4}(a,b) = 0$ satisfying $|\lambda+a| < \lambda$
and $|\lambda+b| < \lambda$.
So $K_4$ does not even have the bivariate property BC${}_\lambda$.
\qed

We can deduce from Lemma~\ref{lemma5.4} and Proposition~\ref{prop5.5}
a necessary and sufficient condition for $G$ to have
various forms of the multivariate Brown--Colbourn property:

\begin{theorem}
  \label{thm5.6}
Let $G$ be a loopless connected graph.  Then the following are equivalent:
\begin{itemize}
   \item[(a)]  $G$ has the multivariate property BC${}_1$.
   \item[(b)]  $G$ has the multivariate property BC${}_\lambda$
       for some $\lambda > 0$.
   \item[(c)]  $G$ is series-parallel.
\end{itemize}
\end{theorem}

\proof
(a) $\Longrightarrow$ (b) is trivial.

(b) $\Longrightarrow$ (c):  Let $G$ be a loopless connected graph
that is not series-parallel.
Then $G$ contains a subgraph $H$ that is a subdivision of $K_4$.\footnote{
   The relevant fact is the following
   \cite[Exercise 8.16 and Proposition 1.7.2]{Diestel_97}:
   $G$ is series-parallel $\Longleftrightarrow$
   $G$ has no $K_4$ minor $\Longleftrightarrow$
   $G$ has no $K_4$ topological minor.
   And the latter statement says precisely that
   $G$ contains no subgraph $H$ that is a subdivision of $K_4$.
   See also \cite{Duffin_65,Oxley_86}.
}
Suppose that $H = (K_4)^{\bowtie \bf s}$ with ${\bf s} = (s_1,\ldots,s_6)$,
and define $s = \max(s_1,\ldots,s_6)$.
Now fix any $\lambda > 0$;
then, by Proposition~\ref{prop5.5}
we can find a vector ${\bf v} = (v_1,\ldots,v_6)$
that is a zero of $C_{K_4}({\bf v})$
and satisfies $|\lambda/s + v_i| < \lambda/s$ for $i=1,\ldots,6$.
It then follows that the vector ${\bf v'} = (v'_1,\ldots,v'_6)$
defined by $v'_i = s_i v_i$
satisfies $C_H({\bf v'}) = 0$
and $|\lambda + v'_i| < \lambda$ for $i=1,\ldots,6$.
Therefore $H$ does not have the multivariate property BC${}_\lambda$.
By Lemma~\ref{lemma5.4}, $G$ cannot have this property either.

(c) $\Longrightarrow$ (a):
This is proven in \cite[Remark 3 in Section 4.1]{Sokal_chromatic_bounds},
but for the convenience of the reader we repeat the proof here.
Suppose that $G$ is a loopless connected series-parallel graph;
this means that $G$ can be obtained from a tree
by a finite sequence of series and parallel extensions of edges
(i.e.\ replacing an edge by two edges in series or two edges in parallel).
We will prove that $G$ has the multivariate property BC${}_1$,
by induction on the length of this sequence of series and parallel extensions.
The base case is when $G$ is a tree:
then $C_G({\bf v}) = \prod_{e \in E(G)} v_e$
and $G$ manifestly has the multivariate property BC${}_1$.
Now suppose that $G$ is obtained from a smaller graph $G'$
by replacing an edge $e_*$ of $G'$ by two parallel edges $e_1,e_2$.
Use the parallel reduction formula \reff{eq.parallel}:
since $|1 + v_1| < 1$ and $|1 + v_2| < 1$ imply $|1 + v_*| < 1$,
we deduce that $G$ has the multivariate property BC${}_1$ if $G'$ does.
Suppose, finally, that $G$ is obtained from a smaller graph $G'$
by replacing an edge $e_*$ of $G'$ by two edges $e_1,e_2$ in series.
Use the series reduction formula \reff{eq.series0}/\reff{eq.series}
and the fact that $|1+v| < 1$ is equivalent to $\real(1/v) < -1/2$:
then $\real(1/v_i) < -1/2$ for $i=1,2$ implies that $\real(1/v_*) < -1 < -1/2$,
and moreover the prefactor $v_1 + v_2$ is nonzero;
so we deduce that $G$ has the multivariate property BC${}_1$ if $G'$ does.
\qed

For each graph $G$, let us define $\lambda_\star(G)$
to be the maximum $\lambda$ for which
$G$ has the multivariate property BC${}_\lambda$.
Then Theorem~\ref{thm5.6}
states a surprising (at first sight) dichotomy:
either $\lambda_\star(G) = 0$ [when $G$ is not series-parallel]
or else $\lambda_\star(G) \ge 1$ [when $G$ is series-parallel].

Some series-parallel graphs have $\lambda_\star(G) = 1$ exactly:
for example, the graphs $K_2^{(n)}$ (a pair of vertices connected
by $n$ parallel edges) have $C_{K_2^{(n)}}(v) = (1+v)^n - 1$
and hence even have univariate roots on the circle $|1+v| = 1$.
On the other hand, some series-parallel graphs have $\lambda_\star(G) > 1$:
for example, the cycles $C_n$ have $\lambda_\star(G) = n/2$.
[{\sc Proof:}  We have
\be
   C_{C_n}({\bf v})  \;=\;
   \left( \prod\limits_{i=1}^n v_i \right)
   \left( 1 + \sum\limits_{i=1}^n {1 \over v_i} \right)
   \;,
\ee
which is nonvanishing if $\real(1/v_i) < -1/n$ for all $i$.
But this is equivalent to $|n/2 + v_i| < n/2$.]
It is an interesting open problem to characterize the graphs
that have $\lambda_\star(G) = 1$ or, more ambitiously,
to find a simple graph-theoretic formula for $\lambda_\star(G)$.


\vspace{2cm}

\section*{Acknowledgments}

We wish to thank Jason Brown, Robert Shrock and Dave Wagner
for valuable conversations and correspondence.
We also wish to thank Dario Bini for supplying us the MPSolve package
\cite{Bini_package,Bini-Fiorentino} for finding roots of polynomials.
Finally, we wish to thank an anonymous referee for helpful suggestions.

This research was supported in part by
U.S.\ National Science Foundation grant PHY--0099393
and by a University of Western Australia research grant.

\end{document}